# Mixing and tight polyhedra

Thomas Ward[1],*

*University of East Anglia*

**Abstract:** Actions of $\mathbb{Z}^d$ by automorphisms of compact zero-dimensional groups exhibit a range of mixing behaviour. Schmidt introduced the notion of mixing shapes for these systems, and proved that non-mixing shapes can only arise non-trivially for actions on zero-dimensional groups. Masser has shown that the failure of higher-order mixing is always witnessed by non-mixing shapes. Here we show how valuations can be used to understand the (non-)mixing behaviour of a certain family of examples. The sharpest information arises for systems corresponding to tight polyhedra.

## 1. Introduction

Let $\alpha$ be a $\mathbb{Z}^d$-action by invertible measure-preserving transformations of a probability space $(X, \mathcal{B}, \mu)$. A sequence of vectors $(\mathbf{n}_1^{(j)}, \mathbf{n}_2^{(j)}, \ldots, \mathbf{n}_r^{(j)})_{j \geqslant 1}$ in $(\mathbb{Z}^d)^r$ that are moving apart in the sense that

$$\mathbf{n}_s^{(j)} - \mathbf{n}_t^{(j)} \longrightarrow \infty \text{ as } j \longrightarrow \infty \text{ for any } s \neq t \tag{1}$$

is called *mixing* for $\alpha$ if for any measurable sets $A_1, \ldots, A_r$,

$$\mu\left(\alpha^{-\mathbf{n}_1^{(j)}}(A_1) \cap \cdots \cap \alpha^{-\mathbf{n}_r^{(j)}}(A_r)\right) \longrightarrow \mu(A_1) \cdots \mu(A_r) \text{ as } j \longrightarrow \infty. \tag{2}$$

If (1) guarantees (2), then $\alpha$ is *r-mixing* or *mixing of order r*. Mixing of order 2 is called simply *mixing*. The maximum value of $r$ for which (1) implies (2) is the *order of mixing* $\mathcal{M}(\alpha)$ of $\alpha$ (if there is no maximum then $\alpha$ is *mixing of all orders*, and we write $\mathcal{M}(\alpha) = \infty$).

For single transformations (the case $d = 1$) it is not known if mixing implies mixing of all orders. For $\mathbb{Z}^2$-actions, Ledrappier's example [4] shows that mixing does not imply 3-mixing. Motivated by the way in which Ledrappier's example fails to be 3-mixing, Schmidt introduced the following notion: A finite set $\{\mathbf{n}_1, \ldots, \mathbf{n}_r\}$ of integer vectors is called a *mixing shape* for $\alpha$ if

$$\mu\left(\alpha^{-k\mathbf{n}_1}(A_1) \cap \cdots \cap \alpha^{-k\mathbf{n}_r}(A_r)\right) \longrightarrow \mu(A_1) \cdots \mu(A_r) \text{ as } k \longrightarrow \infty. \tag{3}$$

The maximum value of $r$ for which (3) holds for all shapes of cardinality $r$ is the *shape order of mixing* $\mathcal{S}(\alpha)$. Clearly $\mathcal{M}(\alpha) \leqslant \mathcal{S}(\alpha)$, but in general there are no other relations; the following is shown in [9].

**Lemma 1.** *For any $s$, $1 \leqslant s \leqslant \infty$, there is a measure-preserving $\mathbb{Z}^2$-action with $\mathcal{M}(\alpha) = 1$ and $\mathcal{S}(\alpha) = s$.*

*The author thanks Manfred Einsiedler for discussions leading to this result.
[1]School of Mathematics, University of East Anglia, Norwich, United Kingdom, e-mail: t.ward@uea.ac.uk

For algebraic systems — those in which $X$ is a compact abelian group, $\mu$ is Haar measure, and $\alpha_\mathbf{n}$ is an automorphism of $X$ for each $\mathbf{n} \in \mathbb{Z}^d$ — if all shapes are mixing, then the system is mixing of all orders (see [6], [8]). Whether the quantitative version of this relationship might hold was asked by Schmidt [7, Problem 2.11]: If all shapes with $r$ elements are mixing, is an algebraic dynamical system $r$-mixing? For $r = 2$ this means that the individual elements of an algebraic $\mathbb{Z}^d$-action are mixing transformations if and only if the whole action is mixing, which is proved in [6, Theorem 1.6]. For $d = 2$ and $r = 3$ this was shown in [2]. Finally, Masser proved this in complete generality [5].

**Theorem 2 (Masser).** *For any algebraic dynamical system $(X, \alpha)$ on a zero-dimensional group $X$, $\mathcal{M}(\alpha) = \mathcal{S}(\alpha)$.*

In conjunction with (4) and the algebraic characterization (5), Theorem 2 shows that $\mathcal{M}(\alpha) = \mathcal{S}(\alpha)$ for any algebraic dynamical system $\alpha$.

The problem of determining the exact order of mixing for a given system remains: By [6, Chap. VIII], there is — in principle — an algorithm that works from a presentation of the module defining an algebraic $\mathbb{Z}^d$-action and determines all the non–mixing shapes, which by Masser's result [5] then determines the exact order of mixing. By [2], all possible orders of mixing arise: for any $m \geqslant 1$ and $d \geqslant 2$, there is an algebraic $\mathbb{Z}^d$-action with $\mathcal{M}(\alpha) = m$

Our purpose here is to show how the methods from [2] extend to $d > 2$. This gives sharp information about mixing properties for a distinguished class of examples associated to tight polyhedra.

## 2. Inequalities for order of mixing

By [8], for an algebraic dynamical system $\alpha$ on a connected group,

$$\mathcal{M}(\alpha) > 1 \Rightarrow \mathcal{M}(\alpha) = \infty, \tag{4}$$

so in particular $\mathcal{M}(\alpha) = \mathcal{S}(\alpha)$ in this case. Thus finite order of mixing for mixing systems can only arise on groups that are not connected. Following [6], any algebraic $\mathbb{Z}^d$-action $\alpha$ on a compact abelian group $X$ is associated via duality to a module $M = M_X$ over the ring $R_d = \mathbb{Z}[u_1^{\pm 1}, \ldots, u_d^{\pm 1}]$ (multiplication by $u_i$ is dual to the automorphism $\alpha^{\mathbf{e}_i}$ for $i = 1, \ldots, d$). Conversely, any $R_d$-module $M$ determines an algebraic $\mathbb{Z}^d$-action $\alpha_M$ on the compact abelian group $X_M$. Approximating the indicator functions of the sets appearing in (2) by finite trigonometric polynomials shows that (2) for $\alpha_M$ is equivalent to the property that for any elements $a_1, \ldots, a_r$ of $M$, not all zero,

$$a_1 \mathbf{u}^{\mathbf{n}_1^{(j)}} + \cdots + a_r \mathbf{u}^{\mathbf{n}_r^{(j)}} = 0_M \tag{5}$$

can only hold for finitely many values of $j$, where

$$\mathbf{u}^\mathbf{n} = u_1^{n_1} \cdots u_d^{n_d}$$

is the monomial corresponding to the position $\mathbf{n} \in \mathbb{Z}^d$. This algebraic formulation of mixing may be used to show that (2) holds for $\alpha_M$ if and only if it holds for all the systems $\alpha_{R_d/\mathfrak{p}}$ for prime ideals $\mathfrak{p}$ associated to $M$ (see [8] for example). The group $X_{R_d/\mathfrak{p}}$ is connected if and only if $\mathfrak{p} \cap \mathbb{Z} = \{0\}$, so these two remarks together mean that it is enough to study systems associated to modules of the form $R_d/\mathfrak{p}$ where $\mathfrak{p}$ is a prime ideal containing a rational prime $p$.

The (dramatic) simplifying assumption made here concerns the shape of the prime ideal $\mathfrak{p}$: from now on, we assume that $\mathfrak{p} = \langle p, \tilde{f} \rangle$ for some polynomial $\tilde{f} \in R_d$. The degree to which this assumption is restrictive depends on $d$: For $d = 2$, any mixing system can be reduced to this case. For $d > 2$, the ideal $\mathfrak{p}$ could take the form $\langle p, \tilde{f}_1, \ldots, \tilde{f}_s \rangle$ for any $s = 1, \ldots, d-1$. In the language of [1], our assumption amounts to requiring that the system be of entropy rank $(d-1)$.

Once the prime $p$ is fixed, the systems we study are therefore parameterized by a single polynomial $\tilde{f} \in R_d$ which is only defined modulo $p$. Since $p$ is fixed, we write $R_{d,p} = \mathbb{F}_p[u_1^{\pm 1}, \ldots, u_d^{\pm 1}]$, and think of the defining polynomial as $f \in R_{d,p}$. Thus the dynamical system we study corresponds to the module

$$R_{d,p}/\langle f \rangle \cong R_d/\langle p, \tilde{f} \rangle \tag{6}$$

where $\tilde{f}$ is any element of $R_d$ with $\tilde{f} \equiv f \pmod{p}$ and the isomorphism in (6) is an isomorphism of $R_d$-modules. Write the polynomial $f$ as a finite sum

$$f(\mathbf{u}) = \sum_{\mathbf{n} \in \mathbb{Z}^d} c_{f,\mathbf{n}} \mathbf{u}^{\mathbf{n}}, \quad c_{f,\mathbf{n}} \in \mathbb{F}_p.$$

The *support* of $f$ is the finite set

$$S(f) = \{\mathbf{n} \in \mathbb{Z}^d \mid c_{f,\mathbf{n}} \neq 0\};$$

denote the convex hull of the support by $N(f)$.

Theorem 2 would follow at once if we knew that a non-mixing *sequence* of order $r$ (that is, a witness to the statement that $\mathcal{M}(\alpha) < r$) was somehow forced to be, or to nearly be, a non-mixing *shape* of order $r$ (a witness to the statement that $\mathcal{S}(\alpha) < r$). The full picture is much more complicated, in part because the presence of the Frobenius automorphism of $\mathbb{F}_p$ leads to many families of solutions to the underlying equations – see [5].

Here we show that in a special setting the simple arguments from [2] do indeed force a non-mixing sequence to approximate a non-mixing shape, giving an elementary approach to Theorem 2 for this very special setting.

Let $P$ be a convex polyhedron in $\mathbb{R}^d$. A *parallel redrawing* of $P$ is another polyhedron $Q$ with the property that every edge of $Q$ is parallel to an edge of $P$. Figure 1 shows a parallel redrawing of a pentagon.

**Definition 3.** A convex polyhedron $P$ in $\mathbb{R}^d$ is *tight* if any parallel redrawing of $P$ is homothetic to $P$.

For example, in $\mathbb{R}^2$, the only tight polyhedra are triangles. In $\mathbb{R}^3$ there are infinitely many combinatorially distinct tight convex polyhedra. Among the Platonic solids, the tetrahedron, octahedron and icosahedron are tight, while the dodecahedron and cube are not. Tightness can be studied via the dimension of the space of parallel redrawings of a polyhedron; see papers of Whiteley [10], [11].

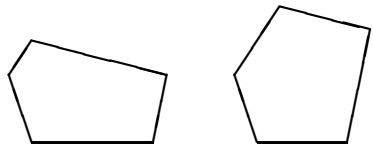

Fig 1. *A parallel redrawing of a pentagon.*

**Theorem 4.** *Let $f$ be an irreducible polynomial in $R_{d,p} = \mathbb{F}_p[u_1^{\pm 1}, \ldots, u_d^{\pm 1}]$, and let $\alpha = \alpha_{R_{d,p}/\langle f \rangle}$ be the algebraic $\mathbb{Z}^d$-action associated to the $R_d$-module $R_{d,p}/\langle f \rangle$. Let $v$ be the number of vertices in $\mathcal{N}(f)$. Then*

(1) *any non-mixing sequence for $\alpha$ along some subsequence contains, with uniform error, a parallel redrawing of $N(f)$;*
(2) *hence $v - 1 \leqslant \mathcal{M}(\alpha) \leqslant \mathcal{S}(\alpha) \leqslant |S(f)| - 1$.*

**Corollary 5.** *If $N(f)$ is tight, then $\mathcal{S}(\alpha) = \mathcal{M}(\alpha)$.*

## 3. Proofs

Throughout we use the characterisation (5) of mixing.

**Lemma 6.** *Let $(\mathbf{n}_1^{(j)}, \ldots, \mathbf{n}_r^{(j)})_{j \geqslant 1}$ be a sequence of $r$-tuples of vectors in $\mathbb{Z}^d$ with the property that there are non-zero elements $a_1, \ldots, a_r \in R_{d,p}/\langle f \rangle$ with*

$$a_1 \mathbf{u}^{\mathbf{n}_1^{(j)}} + \cdots + a_r \mathbf{u}^{\mathbf{n}_r^{(j)}} = 0 \text{ in } R_{d,p}/\langle f \rangle \text{ for all } j \geqslant 1. \qquad (7)$$

*Then there is a constant $K$ with the property that for every edge $e$ of $N(f)$ there is an edge $e'$ of the convex hull of the set $\{\mathbf{n}_1^{(j)}, \ldots, \mathbf{n}_r^{(j)}\}$, joining $\mathbf{n}_s^{(j)}$ to $\mathbf{n}_t^{(j)}$ say, for which there is a point $\widetilde{\mathbf{n}}_s^{(j)}$ with*

(1) $\|\widetilde{\mathbf{n}}_s^{(j)} - \mathbf{n}_s^{(j)}\| \leqslant K$;
(2) *the line through $\widetilde{\mathbf{n}}_s^{(j)} - \mathbf{n}_t^{(j)}$ is parallel to $e$.*

For large $j$, the points $\mathbf{n}_1^{(j)}, \ldots, \mathbf{n}_r^{(j)}$ are widely separated, so Lemma 6 means the edges of the convex hull of these points approximate in direction the edges of $N(f)$ more and more accurately as $j$ goes to infinity.

*Proof of Lemma 6.* Pick an edge $e$ of $N(f)$. Choose a primitive integer vector $\mathbf{v}_1$ orthogonal to $e$ which points outward from $N(f)$ (that is, with the property that for any points $\mathbf{x} \in N(f)$ and $\mathbf{y} \in e$, the scalar product $(\mathbf{x} - \mathbf{y}) \cdot \mathbf{v}_1$ is negative). Also choose an ultrametric valuation $|\cdot|_{\mathbf{v}_1}$ on $R_{d,p}/\langle f \rangle$ with the property that the vector

$$(\log |u_1|_{\mathbf{v}_1}, \ldots, \log |u_d|_{\mathbf{v}_1})^t$$

is a vector of unit length parallel to $\mathbf{v}_1$ that also points outward from $N(f)$. This valuation may be found by extending the vector $\mathbf{v}_1$ to a set of primitive integer vectors $\{\mathbf{v}_1, \mathbf{v}_1^{(2)}, \ldots, \mathbf{v}_1^{(d)}\}$ with $\mathbf{v}_1 \cdot \mathbf{v}_1^{(j)} < 0$ for $j \geqslant 2$ that generates $\mathbb{Z}^d$, as illustrated in Figure 2, and then thinking of $R_{d,p}/\langle f \rangle$ as

$$\mathbb{F}_p[\mathbf{u}^{\mathbf{v}_1}][\mathbf{u}^{\mathbf{v}_1^{(2)}}, \ldots, \mathbf{u}^{\mathbf{v}_1^{(d)}}]/\langle f \rangle.$$

Let

$$K_1 = 2 \max_{i=1,\ldots,r} \{|\log |m_i|_{\mathbf{v}_1}|\}.$$

Now for fixed $j \geqslant 1$ choose $t$ with the property that

$$|\mathbf{u}^{\mathbf{n}_t^{(j)}}|_{\mathbf{v}_1} \geqslant |\mathbf{u}^{\mathbf{n}_s^{(j)}}|_{\mathbf{v}_1} \text{ for all } s, 1 \leqslant s \leqslant r.$$

Then the ultrametric inequality for $|\cdot|_{\mathbf{v}_1}$ and the relation (7) show that there must be (at least) one other vertex $\mathbf{n}_s^{(j)}$ which is no further than $K_1$ from the hyperplane orthogonal to $\mathbf{v}_1$ through $\mathbf{n}_t^{(j)}$.

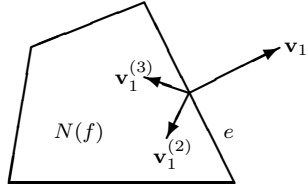

FIG 2. *Extending $\mathbf{v}_1$ to a basis.*

Now choose finitely many vectors $\mathbf{v}_2, \ldots, \mathbf{v}_k$ and a constant $K < \infty$ (depending on the choice of the vectors) with the following property. For each $\ell$, $2 \leqslant \ell \leqslant k$, repeat the construction above corresponding to $\mathbf{v}_1$ and let

$$K_\ell = 2 \max_{i=1,\ldots,r} \{|\log |m_i|_{\mathbf{v}_\ell}|\}.$$

The (purely geometrical) property sought is that any vector $\mathbf{k} \in \mathbb{Z}^d$ with the property that $\mathbf{k}$ is no further than distance $K_\ell$ from the hyperplane orthogonal to $\mathbf{v}_\ell$ through $\mathbf{k}'$ for all $\ell$, $1 \leqslant \ell \leqslant k$, must be within distance $K$ of $\mathbf{k}'$.

Now apply the $k$ different ultrametrics $|\cdot|_{\mathbf{v}_1}, \ldots, |\cdot|_{\mathbf{v}_k}$ to the relation (7) to deduce that there must be a pair of vertices $\mathbf{n}_s^{(j)}$ and $\mathbf{n}_t^{(j)}$ (the parameter $j$ is still fixed; all other quantities including $s$ and $t$ depend on it) with the property that $\mathbf{n}_s^{(j)}$ lies within distance $K_\ell$ of the hyperplane orthogonal to $\mathbf{v}_\ell$ through $\mathbf{n}_t^{(j)}$ for $1 \leqslant \ell \leqslant k$. Since all the vectors $\mathbf{v}_\ell$ are orthogonal to the edge $e$, this proves the lemma. $\square$

*Proof of Theorem 4.* Let $(\mathbf{n}_1^{(j)}, \ldots, \mathbf{n}_r^{(j)})_{j \geqslant 1}$ be a non-mixing sequence for $\alpha$. Thus by (5) there are non-zero elements $a_1, \ldots, a_r \in R_{d,p}/\langle f \rangle$ with

$$a_1 \mathbf{u}^{\mathbf{n}_1^{(j)}} + \cdots + a_r \mathbf{u}^{\mathbf{n}_r^{(j)}} = 0 \text{ in } R_{d,p}/\langle f \rangle \text{ for all } j \geqslant 1. \tag{8}$$

Pick a vertex $v_1$ of $N(f)$ and an edge $e_1$ starting at $v_1$ of $N(f)$, and relabel the non-mixing sequence so that the edge $e_1$ is approximated in direction (in the sense of Lemma 6) by the pair $\mathbf{n}_1^{(j)}, \mathbf{n}_2^{(j)}$ for all $j \geqslant 1$. By Lemma 6, for each $j$ there is a vector $\mathbf{m}_1^{(j)}$ with $\|\mathbf{n}_2^{(j)} - \mathbf{m}_1^{(j)}\| \leqslant K$ such that the line joining $\mathbf{n}_1^{(j)}$ to $\mathbf{m}_1^{(j)}$ is parallel to $e_1$ for all $j$. Since the set of integer vectors $\mathbf{v}$ with $\|\mathbf{v}\| \leqslant K$ is finite, we may find an infinite set $S_1 \subset \mathbb{N}$ with

$$\mathbf{n}_1^{(j)} - \mathbf{m}_1^{(j)} = \mathbf{k}_1, \text{ a constant, for all } j \in S_1.$$

This gives an improved version of the relation (8),

$$a_1 \mathbf{u}^{\mathbf{n}_1^{(j)}} + a_2' \mathbf{u}^{\mathbf{m}_1^{(j)}} + \cdots + a_r \mathbf{u}^{\mathbf{n}_r^{(j)}} = 0 \text{ in } R_{d,p}/\langle f \rangle \text{ for all } j \in S_1 \tag{9}$$

where $a_2' = a_2 \mathbf{u}^{\mathbf{k}_1}$.

Now select another edge $e_2$ of $N(f)$ starting at $v_1$ whose approximating pair is $\mathbf{n}_1^{(j)}$ and (after relabelling) $\mathbf{n}_3^{(j)}$. We now need to allow $2K$ of movement in $\mathbf{n}_3$ to give $\mathbf{m}_3$. This gives an infinite set $S_2 \subset S_1 \subset \mathbb{N}$ and a modified version of (9)

$$a_1 \mathbf{u}^{\mathbf{n}_1^{(j)}} + a_2' \mathbf{u}^{\mathbf{m}_1^{(j)}} + a_3' \mathbf{u}^{\mathbf{m}_2^{(j)}} + \cdots + a_r \mathbf{u}^{\mathbf{n}_r^{(j)}} = 0 \text{ in } R_{d,p}/\langle f \rangle \text{ for all } j \in S_2 \tag{10}$$

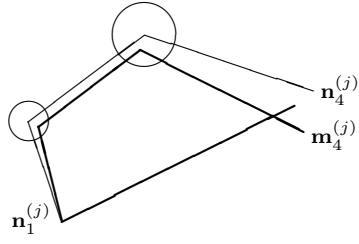

Fig 3. *Approximating a loop in a parallel redrawing of $N(f)$.*

in which $\mathbf{n}_1^{(j)} - \mathbf{m}_1^{(j)}$ is parallel to $e_1$ and $\mathbf{m}_1^{(j)} - \mathbf{m}_2^{(j)}$ is parallel to $e_2$. Continue this process of relabelling, passing to a subsequence and adjusting the coefficients in (10) to exhaust all the edges along some path from $v_1$. The type of situation that may emerge is shown in Figure 3, where $\mathbf{n}_1^{(j)}$ is fixed, $\mathbf{n}_2^{(j)}$ has been moved no further than $K$, $\mathbf{n}_3^{(j)}$ a distance no more than $2K$ and $\mathbf{n}_4^{(j)}$ a distance no more than $3K$ to give edges parallel to edges of $N(f)$. By Lemma 6 there may be an edge of $N(f)$ for which $\mathbf{n}_4^{(j)}$ is the approximating partner, and we have already chosen to adjust $\mathbf{n}_4^{(j)}$ to $\mathbf{m}_4^{(j)}$.

It is difficult to control what loops may arise: for example the Herschel graph [3] shows that a convex polyhedron need not be Hamiltonian as a graph. Nonetheless, the bold path in Figure 3 is, to within a uniformly bounded error, a parallel redrawing of that loop in $N(f)$. This process may be continued to modify all the points $\mathbf{n}_s^{(j)}$ by uniformly bounded amounts to end up with an infinite set $S_* \subset \mathbb{N}$ and a relation

$$a_1 \mathbf{u}^{\mathbf{n}_1^{(j)}} + a_2' \mathbf{u}^{\mathbf{m}_1^{(j)}} + a_3' \mathbf{u}^{\mathbf{m}_2^{(j)}} + \cdots + a_r' \mathbf{u}^{\mathbf{m}_{r-1}^{(j)}} = 0 \text{ in } R_{d,p}/\langle f \rangle \text{ for all } j \in S_* \quad (11)$$

with the property that every edge of $N(f)$ is parallel to within a uniformly bounded error to an edge in the convex hull of the set $\{\mathbf{n}_1^{(j)}, \mathbf{m}_1^{(j)}, \ldots, \mathbf{m}_{r-1}^{(j)}\}$ for all $j \in S_*$, proving part 1. In particular, $r \geqslant v$, so $\mathcal{M}(\alpha) < r$ implies $r \geqslant v$, hence $\mathcal{M}(\alpha) \geqslant v-1$. This proves one of the inequalities in part 2.

All that remains is to prove the other inequality in part 2. If

$$f(\mathbf{u}) = \sum_{\mathbf{n} \in S(f)} c_{f,\mathbf{n}} \mathbf{u}^{\mathbf{n}}, \quad c_{f,\mathbf{n}} \in \mathbb{F}_p$$

then the relation

$$\sum_{\mathbf{n} \in S(f)} c_{f,\mathbf{n}} \mathbf{u}^{\mathbf{n}} = 0 \text{ in } R_{d,p}/\langle f \rangle$$

implies that

$$\Big( \sum_{\mathbf{n} \in S(f)} c_{f,\mathbf{n}} \mathbf{u}^{\mathbf{n}} \Big)^{p^k} = \sum_{\mathbf{n} \in p^k \cdot S(f)} c_{f,\mathbf{n}} \mathbf{u}^{\mathbf{n}} = 0 \text{ in } R_{d,p}/\langle f \rangle \text{ for all } k \geqslant 1,$$

so $S(f)$ is a non-mixing shape, and $\mathcal{S}(\alpha) \leqslant |S(f)| - 1$. $\square$

*Proof of Corollary 5.* If $N(f)$ is tight, then (11) may be improved further: multiply each of the coefficients by a monomial chosen to shift the vertices by a uniformly bounded amount to lie on an integer multiple of $N(f)$. The resulting sequence $\{\tilde{\mathbf{n}}_1^{(j)}, \tilde{\mathbf{m}}_1^{(j)}, \ldots, \tilde{\mathbf{m}}_{r-1}^{(j)}\}$ is homothetic to $N(f)$ and so is a non-mixing shape. Thus $\mathcal{M}(\alpha) < r$ implies that $\mathcal{S}(\alpha) < r$, so $\mathcal{S}(\alpha) \leqslant \mathcal{M}(\alpha)$. $\square$